\documentclass[11pt]{article}

\title{\textbf{
The modular properties and \\
the integral representations of \\
the multiple elliptic gamma functions}}
\author{{\Large Atsushi Narukawa} \\[.1in]
\textit{
Department of Mathematics, School of Science and Engineering,} \\
\textit{
Waseda University, Tokyo 169-8555, Japan}
\footnote{
The author recently transferred to
The Dai-ichi Mutual Life Insurance Company,
Tokyo 100-8411, Japan.} \\[.1in]
\texttt{narukawa@ruri.waseda.jp}}
\date{June 10, 2003}

\usepackage{amssymb}

\newcommand{\dsum}{\displaystyle\sum}
\newcommand{\dfrac}{\displaystyle\frac}

\newcommand{\real}{\mathbb{R}}
\newcommand{\comp}{\mathbb{C}}
\newcommand{\zahl}{\mathbb{Z}}
\newcommand{\im}{\mathrm{Im}\,}
\newcommand{\re}{\mathrm{Re}\,}
\newcommand{\li}{\mathrm{Li}\,}
\newcommand{\res}{\mathop{\mathrm{Res}}}
\newcommand{\widecheck}[1]{\stackrel{\!\!\!\vee}{#1}}
\newcommand{\qed}{\relax\ifmmode\hskip2em
\Box\else\unskip\nobreak\hskip1em $\Box$\fi}

\newtheorem{thm}{Theorem}
\newtheorem{prop}[thm]{Proposition}
\newtheorem{lemma}[thm]{Lemma}
\newtheorem{cor}[thm]{Corollary}

\begin{document}
\maketitle

\begin{abstract}

We show the modular properties of the multiple ``elliptic'' gamma
functions, which are an extension of those of the theta function
and the elliptic gamma function.
The modular property of the theta function is known
as Jacobi's transformation,
and that of the elliptic gamma function was provided
by Felder and Varchenko.
In this paper, we deal with the multiple sine functions,
since the modular properties of the multiple elliptic gamma
functions result from the equivalence between two ways
to represent the multiple sine functions as infinite products.

We also derive integral representations of the multiple sine
functions and the multiple elliptic gamma functions.
We introduce correspondences between the multiple elliptic
gamma functions and the multiple sine functions.

\end{abstract}


\section{Introduction}

The theta function $\theta_0 (z,\tau)$ and
the elliptic gamma function $\Gamma (z,\tau,\sigma)$ are defined
by infinite products
\begin{eqnarray*}
\theta_0 (z,\tau)
&=& \prod_{j=0}^{\infty}
(1-e^{2 \pi i ((j+1)\tau -z)})
(1-e^{2 \pi i (j\tau +z)}) , \\
\Gamma (z,\tau,\sigma)
&=& \prod_{j,k=0}^{\infty}\frac{
1-e^{2 \pi i ((j+1)\tau +(k+1)\sigma -z)}
}{
1-e^{2 \pi i (j\tau +k\sigma +z)}
},
\end{eqnarray*}
then $\Gamma (z,\tau,\sigma)$ satisfies the difference equation
$\Gamma (z+\tau,\tau,\sigma)
= \theta_0 (z,\sigma) \Gamma (z,\tau,\sigma)$.

The elliptic gamma function was originally constructed
by Ruijsenaars \cite{R1} as a unique solution
of difference equations which include the theta function as above.
After his work, Felder and Varchenko \cite{FV} derived
the ``modular property'' of this function,
\begin{equation}
\Gamma \left( \frac{z}{\sigma},
\frac{\tau}{\sigma},-\frac{1}{\sigma} \right)
= e^{ \pi i Q(z;\tau,\sigma) }
\Gamma \left( \frac{z-\sigma}{\tau},
-\frac{\sigma}{\tau},-\frac{1}{\tau} \right)
\Gamma (z,\tau,\sigma),
\label{1}\end{equation}
where
\begin{eqnarray*}
Q (z;\tau,\sigma)
&=& \frac{z^3}{3\tau\sigma}
- \frac{\tau+\sigma-1}{2\tau\sigma} z^2
+ \frac{\tau^2+\sigma^2+1+3\tau\sigma-3\tau-3\sigma}{6\tau\sigma}
 z \\
&& - \frac{(\tau+\sigma-1)(\tau\sigma-\tau-\sigma)}{12\tau\sigma}.
\end{eqnarray*}
This formula is an extension of Jacobi's transformation,
in which the group $\mathrm{SL} (2,\zahl) \ltimes \zahl^2$ acts
on the parameter of the theta function.
Felder and Varchenko deduced (\ref{1})
from the modular properties of $\theta_0 (z,\tau)$
and the special value of $\Gamma (z,\tau,\sigma)$.
They also gave a cohomological interpretation to this formula
with $\mathrm{SL} (3,\zahl) \ltimes \zahl^3$.

On the other hand, Nishizawa \cite{N1} constructed a hierarchy
of meromorphic functions
which includes the theta function and the elliptic gamma function.
He call these new functions the multiple elliptic gamma functions
$G_r (z|\tau_0,\cdots,\tau_r)$, which are considered
as an elliptic analogue of the multiple gamma functions.
They are defined by certain infinite products
called $q$-shifted factorials $(x;\underline{q})^{(r)}_\infty$.
They satisfy functional relations, such as
\[
G_r (z+\tau_j|\tau_0,\cdots,\tau_r)
= G_{r-1} (z|\tau_0,\cdots,\widecheck{\tau_j},\cdots,\tau_r) 
\  G_r (z|\tau_0,\cdots,\tau_r).
\]
Conversely he characterized these functions
with above relations and initial values.

Our main purpose of this paper is to derive the modular properties
of the multiple elliptic gamma functions
$G_r (z|\tau_0,\cdots,\tau_r)$
while we discuss the properties of the multiple sine functions
$S_r (z |\omega_1,\cdots,\omega_r)$
and the multiple Bernoulli polynomials
$B_{r,n} (z |\omega_1,\cdots,\omega_r)$.
The hierarchy of the multiple sine functions is defined
by Barnes' multiple gamma functions.
Those have been studied by Shintani \cite{S},
Kurokawa \cite{K1,K2}.
The multiple Bernoulli polynomials are attached
to the multiple zeta functions and the multiple gamma functions
as in \cite{B,R2}.

We introduce integral representations and
infinite product representations of the multiple sine functions.
Then it is shown that there are two ways
to represent them as an infinite product.
From the equivalence between them, the modular properties of
$G_r (z|\tau_0,\cdots,\tau_r)$ are obtained.
For example, (\ref{1}) is derived from the two representations of
$S_3 (z|\omega_1,\omega_2,\omega_3)$, namely
\begin{eqnarray*}
\lefteqn{S_3 (z|\omega_1,\omega_2,\omega_3)} \\
&=& \exp \left\{
- \frac{\pi i}{6} B_{33} (z|\omega_1,\omega_2,\omega_3)
\right\} \\
&& \times \prod_{j,k=0}^{\infty} \frac{
( 1-e^{2 \pi i
(z/\omega_1 -(j+1) \omega_2/\omega_1 -(k+1) \omega_3/\omega_1) } )
( 1-e^{2 \pi i
(z/\omega_3 +j \omega_1/\omega_3 +k \omega_2/\omega_3) } )
}{
1-e^{2 \pi i
(z/\omega_2 +j \omega_1/\omega_2 -(k+1) \omega_3/\omega_2) }
} \\
&=& \exp \left\{
+ \frac{\pi i}{6} B_{33} (z|\omega_1,\omega_2,\omega_3)
\right\} \\
&& \times \prod_{j,k=0}^{\infty} \frac{
( 1-e^{2 \pi i
(-z/\omega_1 -j \omega_2/\omega_1 -k \omega_3/\omega_1) } )
( 1-e^{2 \pi i
(-z/\omega_3 +(j+1) \omega_1/\omega_3 +(k+1) \omega_2/\omega_3) })
}{
1-e^{2 \pi i
(-z/\omega_2 +(j+1) \omega_1/\omega_2 -k \omega_3/\omega_2) }
}.
\end{eqnarray*}
Subustituting $\omega_1 = \tau ,\omega_2 = \sigma ,\omega_3 = -1$,
we have (\ref{1}) and the fact
\[
Q(z;\tau,\sigma)
= -\frac{1}{3} B_{33} (z|\tau,\sigma,-1).
\]
In general, we use $q$-shifted factorials
$(x;\underline{q})^{(r)}_\infty$ to describe
$G_r (z|\tau_0,\cdots,\tau_r)$
and $S_r (z |\omega_1,\cdots,\omega_r)$.
The $q$-shifted factorial is the exponential of
the generalized $q$-polylogarithm.
The general result which we prove is the the following theorem.

\paragraph*{Theorem}
If $r \ge 2, \im \frac{\omega_j}{\omega_k} \ne 0$,
then the multiple elliptic gamma function satisfies the identity
\[
\prod_{k=1}^r G_{r-2}
\left( \frac{z}{\omega_k} \bigg| \left(
\frac{\omega_1}{\omega_k}, \cdots,
\widecheck{\frac{\omega_k}{\omega_k}}, \cdots,
\frac{\omega_r}{\omega_k} \right) \right)
= \exp \left\{
- \frac{2 \pi i}{r!} B_{rr} (z|\underline{\omega})
\right\}.
\]
\bigskip

The remaining part of this paper is devoted to investigate
the integral representations of $G_r (z|\tau_0,\cdots,\tau_r)$.
We recall the results in \cite{FV,N1} again,
and regard $G_r (z|\tau_0,\cdots,\tau_r)$
as infinite products of
$S_{r+1} (z |\omega_1,\cdots,\omega_{r+1})$.

The paper is organized as follows:
In Section 2 and 3, we review the definition and the properties
of Nishizawa's multiple elliptic gamma functions and those of
the multiple Bernoulli polynomials.
Then in Section 4, we introduce the integral representations
and the infinite product representations of the multiple sine
functions.
In Section 5, we prove the modular properties
of the multiple elliptic gamma functions.
In Section 6, we introduce the integral representations
of the multiple elliptic gamma functions,
and these integrals show
that the multiple elliptic gamma functions are described
as infinite products of the multiple sine functions.


\section{The multiple elliptic gamma functions
$G_r (z | \underline{\tau})$}

In this section, we review the multiple elliptic gamma functions
$G_r (z | \underline{\tau})$ according to Nishizawa \cite{N1}.

Let $x=e^{2 \pi iz}, q_j=e^{2 \pi i \tau_j}$ for
$z \in \comp$ and $\tau_j \in \comp-\real \ (0 \le j \le r)$, and
\begin{eqnarray*}
\underline{q} &=& ( q_0, \cdots\qquad\cdots, q_r ), \\
\underline{q}^- (j)
&=& ( q_0, \cdots, \widecheck{q_j}, \cdots, q_r), \\
\underline{q} [j]
&=& ( q_0, \cdots, q_j^{-1}, \cdots, q_r), \\
\underline{q}^{-1} &=& ( q_0^{-1}, \cdots\qquad\cdots, q_r^{-1} )
, \end{eqnarray*}
where $\widecheck{q_j}$ means the excluding of $q_j$.
When $\im\tau_j >0$ for all $j$, define the $q$-shifted factorial
\[
(x;\underline{q})^{(r)}_\infty
= \prod_{j_0,\cdots,j_r=0}^{\infty}
(1- x q_0^{j_0} \cdots q_r^{j_r}).
\]
This infinite product converges absolutely
when $|q_j| <1$.
Thus this function is a holomorphic function
with regard to $z$, whose zeros are
\[
z = \tau_0 \zahl_{\le 0} +\cdots +\tau_r \zahl_{\le 0} +\zahl.
\]

In general we can define the $q$-shifted factorial
for $\tau_j \in \comp-\real$ as follows:
When $\im\tau_0,\cdots,\im\tau_{k-1} <0$ and
$\im\tau_k,\cdots,\im\tau_r >0$,
that is, $|q_0|,\cdots,|q_{k-1}| > 1$ and
$|q_k|,\cdots,|q_r| < 1$, we difine
\begin{eqnarray}
(x;\underline{q})^{(r)}_\infty
&=& \left\{ (q_0^{-1} \cdots q_{k-1}^{-1} x;
( q_0^{-1}, \cdots, q_{k-1}^{-1}, q_k, \cdots, q_r)
)^{(r)}_\infty \right\}^{(-1)^k} \label{9} \\
&=& \left\{ \prod_{j_0,\cdots,j_r=0}^{\infty}
(1- x q_0^{-j_0 -1} \cdots q_{k-1}^{-j_{k-1} -1}
 q_k^{j_k} \cdots q_r^{j_r}) \right\}^{(-1)^k}. \label{10}
\end{eqnarray}
More general definition can be done in a similar way
as the values are not changed under the permutation of
$q_0, \cdots, q_r$.
In this definition $(x;\underline{q})^{(r)}_\infty$ is
a meromorphic function of $z$ satisfying the following functional
equations.

\begin{prop}\label{1.1}
\[
(x;\underline{q})^{(r)}_\infty
= \frac{1}{(q_j^{-1} x;\underline{q} [j])^{(r)}_\infty}
, \qquad
(q_j x;\underline{q})^{(r)}_\infty
= \frac{(x;\underline{q})^{(r)}_\infty}
{(x;\underline{q}^- (j))^{(r-1)}_\infty}.
\]
\end{prop}

We next denote
\begin{eqnarray*}
\underline{\tau} &=& ( \tau_0, \cdots\qquad\cdots, \tau_r ), \\
\underline{\tau}^- (j)
&=& ( \tau_0, \cdots, \widecheck{\tau_j}, \cdots, \tau_r ), \\
\underline{\tau} [j]
&=& ( \tau_0, \cdots, -\tau_j, \cdots, \tau_r ), \\
- \underline{\tau}
&=& ( -\tau_0, \cdots\qquad\cdots, -\tau_r ), \\
|\underline{\tau}| &=& \tau_0 + \cdots + \tau_r
\end{eqnarray*}
and define the multiple elliptic gamma function
\begin{eqnarray}
G_r (z|\underline{\tau})
&=& (x^{-1} q_0 \cdots q_r ;\underline{q})^{(r)}_\infty
\{(x;\underline{q})^{(r)}_\infty \}^{(-1)^r} \label{13b} \\
&=& \{(x^{-1} ;\underline{q}^{-1})^{(r)}_\infty \}^{(-1)^{r+1}}
\{(x;\underline{q})^{(r)}_\infty \}^{(-1)^r}
\label{13}.
\end{eqnarray}
$G_r (z|\underline{\tau})$ is defined for
$\tau_j \in \comp - \real$ from the general definition of
$(x;\underline{q})^{(r)}_\infty$.
The hierarchy of $G_r (z|\underline{\tau})$ includes
the theta function $\theta_0 (z,\tau) $ and
the elliptic gamma function $\Gamma (z,\tau,\sigma)$
which appeared in \cite{R1,FV}.
When $\im\tau,\im\sigma >0$, recall the definition
\begin{eqnarray*}
\theta_0 (z,\tau)
&=& \prod_{j=0}^{\infty}
(1-e^{2 \pi i ((j+1)\tau -z)})
(1-e^{2 \pi i (j\tau +z)})
= G_0 (z|\tau), \label{2} \\
\Gamma (z,\tau,\sigma)
&=& \prod_{j,k=0}^{\infty}\frac{
1-e^{2 \pi i ((j+1)\tau +(k+1)\sigma -z)}
}{
1-e^{2 \pi i (j\tau +k\sigma +z)}
}
= G_1 (z|\tau,\sigma). \label{5}
\end{eqnarray*}

The definition of $(x;\underline{q})^{(r)}_\infty$ and
Proposition \ref{1.1} imply the functional equations:
\begin{eqnarray}
G_r (z+1|\underline{\tau}) &=& G_r(z|\underline{\tau}),
 \label{14} \\
G_r (z+\tau_j|\underline{\tau})
&=& G_{r-1} (z|\underline{\tau}^- (j)) \ 
G_r (z|\underline{\tau}), \label{15} \\
G_r (z|\underline{\tau})
&=& \frac{1}{G_r (z-\tau_j|\underline{\tau} [j])}, \label{16} \\
G_r (-z|-\underline{\tau})
&=& \frac{1}{G_r (z|\underline{\tau})}, \label{18} \\
G_r (z|\underline{\tau}) G_r (z|\underline{\tau}[j])
&=& \frac{1}{G_{r-1} (z|\underline{\tau}^- (j))}. \label{19}
\end{eqnarray}
$G_r (z|\underline{\tau})$ can be expressed as an infinite product
directly by (\ref{16}) and (\ref{10})
for any $\tau_j \in \comp - \real$.

The zeros and poles are easily observed if $\im\tau_j >0$
for all $j$.
When $r$ is even, $G_r (z|\underline{\tau})$ is holomorphic
on $\comp$, and zeros are written as follows:
\begin{eqnarray*}
zeros &&
z = \tau_0 \zahl_{\le 0} + \cdots + \tau_r \zahl_{\le 0} +\zahl,
 \\ &&
z = \tau_0 \zahl_{\ge 1} + \cdots + \tau_r \zahl_{\ge 1} +\zahl.
\end{eqnarray*}
When $r$ is odd, $G_r (z|\underline{\tau})$ is meromorphic
on $\comp$ with poles and zeros written as follows:
\begin{eqnarray*}
poles &
z = \tau_0 \zahl_{\le 0} + \cdots + \tau_r \zahl_{\le 0} +\zahl,
 \\
zeros &
z = \tau_0 \zahl_{\ge 1} + \cdots + \tau_r \zahl_{\ge 1} +\zahl.
\end{eqnarray*}
In particular, $G_r (z|\underline{\tau})$ has
no poles and no zeros in the domain
$\{ 0< \im z < \im |\underline{\tau}| \}$
if $\im\tau_j >0$ for all $j$.

\bigskip

When $r=0$, $G_0 (z|\tau_0)$ means the theta function
$\theta_0 (z,\tau)$.
As we know, $\theta_0 (z,\tau)$ posesses the periodicity
$\theta_0 (z+1,\tau) = \theta_0 (z,\tau),\ 
\theta_0 (z+\tau,\tau)
= e^{ -2 \pi i ( z-1/2 ) } \theta_0 (z,\tau) $
and the modular property
\begin{equation}
\theta_0 \left( \frac{z}{\tau},-\frac{1}{\tau} \right)
= \exp \left\{ \pi i \left(
\frac{z^2}{\tau} + \frac{z}{\tau} -z
+\frac{\tau}{6} + \frac{1}{6\tau} - \frac{1}{2}
\right) \right\} \theta_0 (z,\tau).
\label{3}\end{equation}
The modular property of $r=1$ case has been already described
in (\ref{1}) refering to \cite{FV}.
We derive the modular properties for general
$G_r (z|\underline{\tau})$ in this paper.


\section{The multiple Bernoulli polynomials
$B_{r,n} (z | \underline{\omega})$}

For $z \in \comp$,
\ $\underline{\omega}=(\omega_1,\cdots,\omega_r)$,
\ $\omega_j \in \comp - \{ 0 \}$,
we define the multiple Bernoulli polynomials
$B_{r,n} (z | \underline{\omega})$ with a generating function
\[
\frac{t^r e^{zt}}{\prod_{j=1}^{r} (e^{\omega_j t} -1)}
= \sum_{n=0}^{\infty}
B_{r,n} (z | \underline{\omega}) \frac{t^n}{n!}.
\]
They essentially appeared in \cite{B} and play an important role
to study the multiple zeta functions and
the multiple gamma functions as in the next section.

$B_{r,n} (z | \underline{\omega})$ is a polynomial of degree
$n$ in $z$ and is symmetric in $\omega_1,\cdots,\omega_r$.
It is easy to show that
\begin{eqnarray}
B_{r,n} (cz | c \underline{\omega})
&=& c^{n-r} B_{r,n} (z | \underline{\omega})
\qquad (\forall c \in \comp - \{ 0 \} ), \label{20} \\
B_{r,n} (|\underline{\omega}| -z | \underline{\omega})
&=& (-1)^n B_{r,n} (z | \underline{\omega}), \label{21} \\
B_{r,n} (z+\omega_j | \underline{\omega})
-B_{r,n} (z | \underline{\omega})
&=& n B_{r-1,n-1} (z | \underline{\omega}^{-}(j)), \label{22} \\
B_{r,n} (z | \underline{\omega}[j])
&=& - B_{r,n} (z + \omega_j | \underline{\omega}), \label{23} \\
B_{r,n} (z | \underline{\omega})
+B_{r,n} (z | \underline{\omega}[j])
&=& -n B_{r-1,n-1} (z | \underline{\omega}^{-}(j)), \label{29} \\
\frac{d}{dz} B_{r,n} (z | \underline{\omega})
&=& n B_{r,n-1} (z | \underline{\omega})
\end{eqnarray}
where
\begin{eqnarray*}
c \underline{\omega}
&=& (c \omega_1,\cdots, c \omega_r), \\
|\underline{\omega}|
&=& \omega_1 + \cdots + \omega_r, \\
\underline{\omega}^{-}(j)
&=& (\omega_1,\cdots,\widecheck{\omega_{j}},\cdots,\omega_r)
, \\
\underline{\omega}[j]
&=& (\omega_1,\cdots,-\omega_j,\cdots,\omega_r)
\end{eqnarray*}
and $\widecheck{\omega_{j}}$ means the excluding of
$\omega_{j}$.
In particular, $B_{r,n} (z) = B_{r,n} (z | 1,\cdots,1 )$
obeys
\[
r B_{r+1,n} (z+1)
= (r-n) B_{r,n} (z)
+ nz B_{r,n-1} (z),
\]
\[
B_{r+1,r} (z)
= (z-1) \cdots (z-r)
= r! {z-1 \choose r}.
\label{25}\]

For example, we can see that \newpage
\begin{eqnarray*}
B_{11} (z|\omega_1)
&=& \frac{z}{\omega_1}-\frac{1}{2}, \label{26} \\
B_{22} (z|\omega_1,\omega_2)
&=& \frac{z^2}{\omega_1 \omega_2}
- \frac{\omega_1 + \omega_2}{\omega_1 \omega_2} z
+ \frac{\omega_1^2 + \omega_2^2 + 3 \omega_1 \omega_2}
{6 \omega_1 \omega_2}, \label{27} \\
B_{33} (z|\omega_1,\omega_2,\omega_3)
&=& \frac{z^3}{\omega_1 \omega_2 \omega_3}
- \frac{3(\omega_1 + \omega_2 + \omega_3)}
{2 \omega_1 \omega_2 \omega_3} z^2 \\
&& + \frac{\omega_1^2 + \omega_2^2 + \omega_3^2
+ 3 \omega_1 \omega_2 + 3 \omega_2 \omega_3 + 3 \omega_3 \omega_1}
{2 \omega_1 \omega_2 \omega_3} z \nonumber \\
&& - \frac{(\omega_1 + \omega_2 + \omega_3)
(\omega_1 \omega_2 + \omega_2 \omega_3 + \omega_3 \omega_1)}
{4 \omega_1 \omega_2 \omega_3}.
\end{eqnarray*}


\section{Definition and properties of the multiple sine functions
$S_r (z | \underline{\omega})$}

\subsection{Definition of the multiple sine functions
$S_r (z | \underline{\omega})$}

Now suppose that the points representing
$\omega_1,\cdots,\omega_r \in \comp$ all lie on the same side of
some straight line through the origin.
(Usually we suppose $\omega_1,\cdots,\omega_r \in \comp$ all lie
on the right half plane.)
In this case, the multiple zeta function is defined by the series
\[
\zeta_r (s,z | \underline{\omega})
= \sum_{n_1,\cdots,n_r=0}^\infty
\frac{1}{(n_1 \omega_1 + \cdots + n_r \omega_r +z)^s}
\]
for $z \in \comp, \re s >r$, where the exponential is rendered
one-valued.
This series is holomorphic in the domain $\{ \re s >r \}$,
and it is analitically continued to $s \in \comp$.
Since it is  holomorphic at $s=0$,
we can next define the multiple gamma function by
\begin{equation}
\Gamma_r (z | \underline{\omega})
= \exp \left( \frac{\partial}{\partial s}
\zeta_r (s,z | \underline{\omega}) \Big|_{s=0} \right).
\label{30}\end{equation}
Now we define the multiple sine function by the form
\begin{equation}
S_r (z | \underline{\omega})
= \Gamma_r (z | \underline{\omega})^{-1}
\Gamma_r (|\underline{\omega}| -z | \underline{\omega})^{(-1)^r}.
\label{31}\end{equation}

The above definition of $\zeta_r (s,z | \underline{\omega})$ is
due to Barnes \cite{B}, and the definitions of
$\Gamma_r (z | \underline{\omega})$ and
$S_r (z | \underline{\omega})$ are due to Kurokawa \cite{K2}
and Jimbo-Miwa \cite{JM}.
Barnes' $\Gamma_r (z | \underline{\omega})$ is slightly
different in the coefficient called Barnes' modular constant.

In particular, the double sine function
$S_2 (z | \omega_1,\omega_2)$ have been studied to construct
solutions or operators of certain equations
of mathematical physics as in \cite{JM} or \cite{KLS}.

\subsection{Integral representations of
$S_r (z|\underline{\omega})$}

\begin{prop}\label{3.1}
{\rm (i)} \ 
$S_r ( c z | c \underline{\omega})
= S_r (z|\underline{\omega})$
for all $c \in \comp - \{ 0 \}$. \\
{\rm (ii)} \ 
When $0 < \re\omega_j \ (\forall j)$ and
$0 < \re z < \re |\underline{\omega}|$,\ 
$S_r (z|\underline{\omega})$ has the integral representations
\begin{eqnarray}
S_r (z|\underline{\omega})
&=& \exp\left\{
(-1)^r \frac{\pi i}{r!} B_{rr} (z|\underline{\omega})
+ (-1)^r \int_{\real +i0}
\frac{e^{zt}}{t \prod_{j=1}^r (e^{\omega_j t}-1)} dt
\right\} \label{32} \\
&=& \exp\left\{
(-1)^{r-1} \frac{\pi i}{r!} B_{rr} (z|\underline{\omega})
+ (-1)^r \int_{\real -i0}
\frac{e^{zt}}{t \prod_{j=1}^r (e^{\omega_j t}-1)} dt
\right\}, \label{33}
\end{eqnarray}
where the contours are taken as following figures.
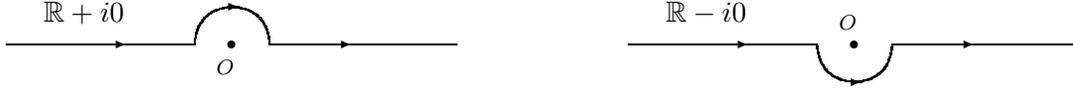
\begin{center}
\setlength{\unitlength}{1mm}
\begin{picture}(60,20)(-30,-10)
\put(0,0){\circle*{1}}
\put(-2,-4){$\scriptstyle O$}
\put(-30,0){\line(1,0){25}}
\put(5,0){\line(1,0){25}}
\bezier{100}(-5,0)(-4.7,4.7)(0,5)
\bezier{100}(5,0)(4.7,4.7)(0,5)
\put(-15,0){\vector(1,0){1}}
\put(15,0){\vector(1,0){1}}
\put(0,5){\vector(1,0){1}}
\put(-25,3){$\real +i0$}
\end{picture}
\hspace{20mm}
\begin{picture}(60,20)(-30,-10)
\put(0,0){\circle*{1}}
\put(-2,2){$\scriptstyle O$}
\put(-30,0){\line(1,0){25}}
\put(5,0){\line(1,0){25}}
\bezier{100}(-5,0)(-4.7,-4.7)(0,-5)
\bezier{100}(5,0)(4.7,-4.7)(0,-5)
\put(-15,0){\vector(1,0){1}}
\put(15,0){\vector(1,0){1}}
\put(0,-5){\vector(1,0){1}}
\put(-25,3){$\real -i0$}
\end{picture}
\end{center}

{\rm (iii)} \ 
When $0 < \re\omega_j \ (\forall j)$,
$S_r (z|\underline{\omega})$ has no poles and no zeros 
in the domain $\{ 0 < \re z < \re |\underline{\omega}| \}$.
\end{prop}

\paragraph*{Proof}

Let $L$ be the half-line from the origin,
let $\overline{L}$ be the half-line conjugate to $L$ with respect
to the real axis,
and let $\overline{L}^\perp$ be the line at right angles
to $\overline{L}$.
Further we take a contour $\tilde{L}$ which embraces $L$,
and other contours $C_0, C$ as following figures.

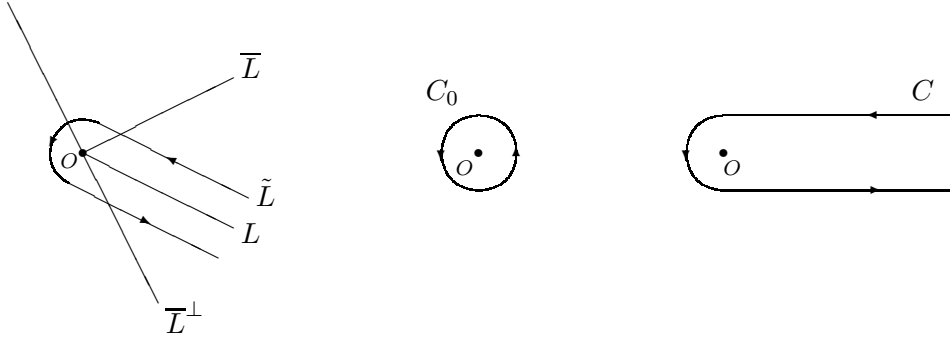
\begin{center}
\setlength{\unitlength}{1mm}
\begin{picture}(40,40)(-20,-20)
\put(0,0){\circle*{1}}
\put(-3,-2){$\scriptstyle O$}
\put(0,0){\line(2,-1){20}}
\put(21,-12){$L$}
\put(2,4){\line(2,-1){20}}
\put(-2,-4){\line(2,-1){20}}
\bezier{100}(2,4)(-2,5.5)(-4,2)
\bezier{100}(-2,-4)(-5.5,-2)(-4,2)
\put(-3.5,2.8){\vector(-1,-2){1}}
\put(12,-1){\vector(-2,1){1}}
\put(8,-9){\vector(2,-1){1}}
\put(23,-7){$\tilde{L}$}
\put(0,0){\line(2,1){20}}
\put(21,10){$\overline{L}$}
\put(0,0){\line(1,-2){10}}
\put(0,0){\line(-1,2){10}}
\put(11,-24){$\overline{L}^\perp$}
\end{picture}
\hspace{20mm}
\begin{picture}(20,40)(-10,-20)
\put(0,0){\circle*{1}}
\put(-3,-3){$\scriptstyle O$}
\put(-7,7){$C_0$}
\bezier{100}(5,0)(4.7,4.7)(0,5)
\bezier{100}(-5,0)(-4.7,4.7)(0,5)
\bezier{100}(-5,0)(-4.7,-4.7)(0,-5)
\bezier{100}(5,0)(4.7,-4.7)(0,-5)
\put(5,0){\vector(0,1){1}}
\put(-5,0){\vector(0,-1){1}}
\end{picture}
\hspace{10mm}
\begin{picture}(40,40)(-10,-20)
\put(0,0){\circle*{1}}
\put(0,-3){$\scriptstyle O$}
\put(0,5){\line(1,0){30}}
\put(0,-5){\line(1,0){30}}
\bezier{100}(-5,0)(-4.7,4.7)(0,5)
\bezier{100}(-5,0)(-4.7,-4.7)(0,-5)
\put(20,5){\vector(-1,0){1}}
\put(20,-5){\vector(1,0){1}}
\put(-5,0){\vector(0,-1){1}}
\put(25,7){$C$}
\end{picture}
\end{center}

Now we assume all of $z,\omega_1,\cdots,\omega_r$ and
$\overline{L}$ lie on the same side with respect to
$\overline{L}^\perp$.
Then as in \cite{B}, $\Gamma_r (z | \underline{\omega})$ has
an integral representation
\[
\Gamma_r (z|\underline{\omega})
= \exp\left\{ \frac{1}{2 \pi i}
\int_{\tilde{L}} \frac{
e^{-zt} \{ \log (-t) + \gamma \}
}{
t \prod_{j=1}^r (1-e^{-\omega_j t})} dt
\right\},
\]
where $\log (-t)$ is rendered one-valued by the cross-cut along $L$, and $\log (-t)$ is real
when $t$ is in $\real_{<0}$ (real and negative).
$\gamma$ is Euler's constant.

Under this assumption, study $S_r ( cz | c \underline{\omega})$
for any $c \in \comp$ with $-c \notin L$.
Then all of $cz$, $c \omega_1$, $\cdots$, $c \omega_r$
lie on the same side with a half-line
$c \overline{L} = \{ ct | \ t \in \overline{L} \}$.
Let $c^{-1} L$ be the half-line conjugate
to $c \overline{L}$, and let $c^{-1} \tilde{L}$ be a contour
which embraces the half-line $c^{-1} L$.
Then we obtain the rotated expression
\[
\Gamma_r (cz | c \underline{\omega})
= \exp\left\{ \frac{1}{2 \pi i}
\int_{c^{-1} \tilde{L}} \frac{
e^{-czt} \{ \log (-t) + \gamma \}
}{
t \prod_{j=1}^r (1-e^{-c \omega_j t})} dt
\right\},
\]
where the cross-cut of $\log (-t)$ is the half-line $c^{-1} L$,
and $\log (-t)$ is real when $t \in \real_{<0}$.
Further changing $t$ into $c^{-1} t$ and changing the branch
of the logarithm, we get
\[
\Gamma_r (cz | c \underline{\omega})
= \exp\left\{ \frac{1}{2 \pi i}
\int_{\tilde{L}} \frac{
e^{-zt} \{ \log (-t) -\log c +\gamma \}
}{
t \prod_{j=1}^r (1-e^{-\omega_j t})} dt
\right\},
\]
where the cross-cut of $\log (-t)$ is the half-line $L$,
and $\log (-t)$ is real when $t \in \real_{<0}$.

If we denote
$\varphi (t)
= \dfrac{e^{zt}}{t \prod_{j=1}^r (e^{\omega_j t}-1)}$,
the above expression is rewritten into
\[
\Gamma_r (cz | c \underline{\omega})
= \exp\left\{\int_{\tilde{L}} (-1)^{r+1} \varphi (-t)
\frac{\{ \log (-t) -\log c +\gamma \}}{2 \pi i}
dt \right\}.
\]
Moreover assuming that $|\underline{\omega}|-z$ lies
on the same side with $\overline{L}$
as well as $z$ and $\underline{\omega}$, we get
\[
\Gamma_r (c |\underline{\omega}| -cz | c \underline{\omega})
= \exp\left\{\int_{\tilde{L}} \varphi (t)
\frac{\{ \log (-t) -\log c +\gamma \}}{2 \pi i}
dt \right\}.
\]
Thus by the definition (\ref{31}),
\[
S_r (cz | c \underline{\omega})
= \exp\left\{ (-1)^r \int_{\tilde{L}}
\{ \varphi (t) + \varphi (-t) \}
\frac{\{ \log (-t) -\log c +\gamma \}}{2 \pi i}
dt \right\}.
\]
When we consider the integral
$\int_{\tilde{L}} \varphi (t) + \varphi (-t) dt$,
$\tilde{L}$ can be replaced by $C_0$
because $\varphi (t)$ and $\varphi (-t)$ are rapidly decreasing
along the half-line $L$.
We note that $\varphi (t) + \varphi (-t)$ is an even function,
then we have
\begin{eqnarray*}
\int_{\tilde{L}} \varphi (t) + \varphi (-t) dt
= \int_{C_0} \varphi (t) + \varphi (-t) dt
=0.
\end{eqnarray*}
They yield
\begin{equation}
S_r (cz | c \underline{\omega})
= \exp\left\{ (-1)^r \int_{\tilde{L}}
\{ \varphi (t) + \varphi (-t) \}
\frac{\log (-t)}{2 \pi i}
dt \right\}. \label{35}
\end{equation}
The right-hand side of this formula is independent of $c$.
This means that $S_r ( c z | c \underline{\omega})$ coincides with
$S_r (z|\underline{\omega})$
under the assumption that $z$ and $|\underline{\omega}|-z$ lie
on the same side with $\overline{L}$.
Nevertheless the analytic continuation certifies that
$S_r ( c z | c \underline{\omega})
= S_r (z|\underline{\omega})$
for any $z$ except the poles.
Another assumption $-c \notin L$ is not essential
since we could choose a line $L$ such that $-c \notin L$
for any $c$.
We conclude (i).

Next we have to show (\ref{32}) of (ii).
Suppose $0 < \re\omega_j$, $0 < \re z < \re |\underline{\omega}|$,
$c=1$, $L= \real_{\ge 0}, \tilde{L} =C$.
We shall transform the integral of (\ref{35}) after fixing
the branch of the logarithm such that
the cross-cut of $\log (-t)$ is the half-line $\real_{\ge 0}$,
and $\log (-t)$ is real when $t \in \real_{<0}$.

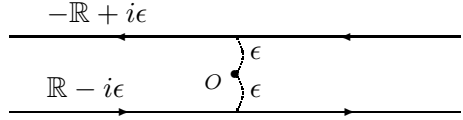
\begin{center}
\setlength{\unitlength}{1mm}
\begin{picture}(60,15)(-30,-5)
\put(0,0){\circle*{1}}
\put(-4,-2){$\scriptstyle O$}
\put(-25,7){$- \real +i \epsilon$}
\put(-30,5){\line(1,0){60}}
\put(-15,5){\vector(-1,0){1}}
\put(15,5){\vector(-1,0){1}}
\bezier{15}(0,5)(2,2.5)(0,0)
\put(2,2){$\epsilon$}
\put(-25,-3){$\real -i \epsilon$}
\put(-30,-5){\line(1,0){60}}
\put(-15,-5){\vector(1,0){1}}
\put(15,-5){\vector(1,0){1}}
\bezier{15}(0,-5)(2,-2.5)(0,0)
\put(2,-3){$\epsilon$}
\end{picture}
\end{center}

Set the contours as above figure,
then stretching the contour $C$ implies
\[
\int_C \varphi (-t) \frac{\log (-t)}{2 \pi i} dt
= \int_{-\real +i \epsilon} \varphi (-t)
 \frac{\log (-t)}{2 \pi i} dt
+\int_{\real -i \epsilon} \varphi (-t)
 \frac{\log (-t)}{2 \pi i} dt.
\]
By changing $t$ into $-t$ without changing the branch
of the logarithm, we have
\begin{eqnarray*}
\int_{-\real +i \epsilon} \varphi (-t)
 \frac{\log (-t)}{2 \pi i} dt
&=& - \int_{\real -i \epsilon} \varphi (t)
 \frac{\log (-t) - \pi i}{2 \pi i} dt, \\
\int_{\real -i \epsilon} \varphi (-t)
 \frac{\log (-t)}{2 \pi i} dt
&=& - \int_{- \real +i \epsilon} \varphi (t)
 \frac{\log (-t) + \pi i}{2 \pi i} dt.
\end{eqnarray*}
Thus
\[
\int_C \varphi (-t) \frac{\log (-t)}{2 \pi i} dt
= - \int_{\real -i \epsilon} \varphi (t)
 \frac{\log (-t) - \pi i}{2 \pi i} dt
- \int_{- \real +i \epsilon} \varphi (t)
 \frac{\log (-t) + \pi i}{2 \pi i} dt.
\]
On the other hand, stretching the contour $C$ implies
\[
\int_C \varphi (t) \frac{\log (-t)}{2 \pi i} dt
= \int_{-\real +i \epsilon} \varphi (t)
 \frac{\log (-t)}{2 \pi i} dt
+\int_{\real -i \epsilon} \varphi (t)
 \frac{\log (-t)}{2 \pi i} dt.
\]
Consequently we obtain
\begin{eqnarray*}
\int_C \{ \varphi (t) + \varphi (-t) \}
\frac{\log (-t)}{2 \pi i} dt
&=& - \frac{1}{2} \int_{-\real +i \epsilon} \varphi (t) dt
+ \frac{1}{2} \int_{\real -i \epsilon} \varphi (t) dt \\
&=& \frac{1}{2} \int_{\real +i0} \varphi (t) dt
+ \frac{1}{2} \left\{ \int_{\real +i0} \varphi (t) dt
+ \int_{C_0} \varphi (t) dt \right\} \\
&=& \int_{\real +i0} \varphi (t) dt
+ \frac{\pi i}{r!} B_{rr} (z | \underline{\omega}),
\end{eqnarray*}
since the residue of $\varphi (t)$ at $t=0$
is $\frac{1}{r!} B_{rr} (z | \underline{\omega})$.
(\ref{32}) is shown.

Eq.(\ref{33}) follows easily from (\ref{32}).
(iii) is due to (ii) because the integrals of (ii) are bounded.
 \qed

\bigskip

The following well-known formulae hold from
the integral representations of $S_r (z | \underline{\omega})$.
\[
S_r (z + \omega_j | \underline{\omega})
= S_{r-1} (z | \underline{\omega}^{-}(j))^{-1}
S_r (z | \underline{\omega}),
\qquad
S_r (z | \underline{\omega})
S_r (|\underline{\omega}|-z | \underline{\omega})^{(-1)^r}
= 1.
\]

\subsection{Definition of the generalized $q$-polylogarithm 
$\li_{r+2} (x;\underline{q})$}

To describe the infinite product representations of
$S_r (z|\underline{\omega})$, we start
with a definition and a lemma.
Now let $x=e^{2 \pi iz}, q_j=e^{2 \pi i \tau_j}$.
When $\im z >0$ and $\im\tau_j \ne 0 \  (0 \le j \le r)$,
Nishizawa \cite{N1} defined the generalized $q$-polylogarithm
\[
\li_{r+2} (x;\underline{q})
= \sum_{n=1}^\infty \frac{x^n}{n \prod_{j=0}^r (1-q_j^n)}.
\]
This series converges absolutely
and then it is holomorphic in $z$.
This function is a generalization
of Kirillov's quantum polylogarithm \cite{Ki},
whose parameters $q_0, \cdots, q_r$ are all equal.
If $r=0$, the above series is
the quantum dilogarithm in \cite{FK}.

We can see the functional equation
\begin{equation}
\li_{r+2} (x;\underline{q})
= - \li_{r+2} (q_j^{-1} x;\underline{q} [j]). \label{38}
\end{equation}
Let us review the following fundamental lemma
verified in \cite{N1}.
\begin{lemma}\label{3.6}
Let $\im z >0, \im\tau_j \ne 0$, that is,
$|x| <1, |q_j| \ne 1$, then
\[
(x;\underline{q})^{(r)}_\infty
= \exp(-\li_{r+2} (x;\underline{q})).
\]
\end{lemma}

\paragraph*{Proof}
First we discuss the case when $|q_j| <1 \ ( \forall j )$.
Applying the formula
$1-x = \exp\left( -\dsum_{n=1}^\infty \frac{x^n}{n} \right) $,
it follows that
\begin{eqnarray*}
\exp(-\li_{r+2} (x;\underline{q}))
&=& \exp\left\{
-\sum_{n=1}^\infty \frac{x^n}{n \prod_{j=0}^r (1-q_j^n)}
\right\} \\
&=& \exp\left\{
- \sum_{n=1}^\infty \sum_{j_0,\cdots,j_r=0}^\infty
\frac{x^n q_0^{n j_0} \cdots q_r^{n j_r}}{n}
\right\} \\
&=& \prod_{j_0,\cdots,j_r=0}^\infty \exp\left\{
- \sum_{n=1}^\infty
\frac{(x q_0^{j_0} \cdots q_r^{j_r} )^n}{n}
\right\} \\
&=& \prod_{j_0,\cdots,j_r=0}^\infty
(1- x q_0^{j_0} \cdots q_r^{j_r} ) \\
&=& (x;\underline{q})^{(r)}_\infty.
\end{eqnarray*}

Next it is enough to discuss the case when
$|q_0|,\cdots,|q_{k-1}| > 1$, and $|q_k|,\cdots,|q_r| < 1$.
The first case and (\ref{9}),(\ref{38}) implies
\begin{eqnarray*}
(x;\underline{q})^{(r)}_\infty
&=& \{ (q_0^{-1} \cdots q_{k-1}^{-1} x;
( q_0^{-1}, \cdots, q_{k-1}^{-1}, q_k, \cdots, q_r)
)^{(r)}_\infty \}^{(-1)^k} \\
&=& \{ \exp(-\li_{r+2} (q_0^{-1} \cdots q_{k-1}^{-1} x;
( q_0^{-1}, \cdots, q_{k-1}^{-1}, q_k, \cdots, q_r)
) \}^{(-1)^k} \\
&=& \{ \exp(- (-1)^k \ \li_{r+2} (x;
( q_0, \cdots, q_{k-1}, q_k, \cdots, q_r)
) \}^{(-1)^k} \\
&=& \exp(-\li_{r+2} (x;\underline{q})).
\qed \end{eqnarray*}

\subsection{Infinite product representations of
$S_r (z|\underline{\omega})$}

We need the following lemma.

\begin{lemma}\label{3.7}
If $ 0< \re z < \re |\underline{\omega}| $,\ 
$ \im z >0 $ and $ \im z > \im |\underline{\omega}| $,
there exists a real series $\{ a_n \}$ such that
\[
\lim_{n \to \infty} a_n = + \infty \quad \mbox{and} \quad
\lim_{n \to + \infty} \int_{\real +i a_n}
\frac{e^{zt}}{t \prod_{j=1}^r (e^{\omega_j t}-1)} dt
=0.
\]

Similarly if $ 0< \re z < \re |\underline{\omega}| $,\ 
$ \im z <0 $ and $ \im z < \im |\underline{\omega}| $,
there exists a real series $\{ a_n \}$ such that
\[
\lim_{n \to \infty} a_n = - \infty \quad \mbox{and} \quad
\lim_{n \to + \infty} \int_{\real +i a_n}
\frac{e^{zt}}{t \prod_{j=1}^r (e^{\omega_j t}-1)} dt
=0,
\]
where the contour is drawn in the following figure.
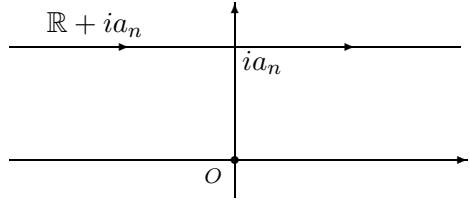
\begin{center}
\setlength{\unitlength}{1mm}
\begin{picture}(60,25)(-30,-5)
\put(0,0){\circle*{1}}
\put(-4,-3){$\scriptstyle O$}
\put(-25,17){$\real +i a_n$}
\put(-30,0){\line(1,0){60}}
\put(-30,15){\line(1,0){60}}
\put(0,-5){\line(0,1){25}}
\put(30,0){\vector(1,0){1}}
\put(0,20){\vector(0,1){1}}
\put(-15,15){\vector(1,0){1}}
\put(15,15){\vector(1,0){1}}
\put(1,12){$i a_n$}
\end{picture}
\end{center}
\end{lemma}

There exist a small $\varepsilon > 0$ and
a real series $\{ a_n \}$ such that
the distances from each poles to any contours $\real +i a_n$
are more than $\varepsilon$.
We reach the above lemma
by estimating the absolute value of the integrand.

Now we set
$x_k=e^{2 \pi iz / \omega_k}$,
$q_{jk}=e^{2 \pi i \omega_j / \omega_k}$,
$\underline{q_k}
= ( q_{1k}, \cdots, \widecheck{q_{kk}}, \cdots, q_{rk} )$
and
$\underline{q_k^{-1}}
= ( q_{1k}^{-1}, \cdots, \widecheck{q_{kk}^{-1}},
 \cdots, q_{rk}^{-1} )$,
then we are ready to prove the formula on
$S_r (z|\underline{\omega})$.
\begin{prop}\label{3.8}
If $r \ge 2, \im \frac{\omega_j}{\omega_k} \ne 0$,
then $S_r (z|\underline{\omega})$ has
the following infinite product representations.
\begin{eqnarray}
S_r (z|\underline{\omega})
&=& \exp \left\{
(-1)^r \frac{\pi i}{r!} B_{rr} (z|\underline{\omega})
\right\}
\prod_{k=1}^{r} (x_k ; \underline{q_k})^{(r-2)}_\infty
 \label{39} \\
&=& \exp \left\{
(-1)^{r-1} \frac{\pi i}{r!} B_{rr} (z|\underline{\omega})
\right\}
\prod_{k=1}^{r} (x_k^{-1} ; \underline{q_k^{-1}})^{(r-2)}_\infty.
 \label{40}
\end{eqnarray}
\end{prop}

\paragraph*{Proof}
From Proposition \ref{3.1}(i) and (\ref{20}),
it is enough to discuss the case $0 < \re\omega_j$ for all $j$.
We need to evaluate the integrals of Proposition \ref{3.1}(ii)
by the residue formula.
We first computate that
\[
\res_{t=\frac{2 \pi in}{\omega_k}}
\frac{e^{zt}}{t \prod_{j=1}^r (e^{\omega_j t}-1)}
= \frac{x_k^n}
{2 \pi in \prod_{j=1, j \ne k}^r (q_{jk}^n-1)}.
\]
For the convergency of the series,
we restrict the domain of $z$ to
\[
\left\{ z \in \comp \Big|
0 < \re z < \re |\underline{\omega}| \ ,\ 
0 < \im z \ ,\ 
\im |\underline{\omega}| < \im z \ ,\ 
\im \frac{z}{\omega_k} >0 \ (\forall k)
\right\},
\]
which is not void when $\re \omega_j >0$.
Then we can use Lemma \ref{3.7},
and add up the residues in $\{ \im t >0 \}$, namely
\begin{eqnarray*}
\lefteqn{\exp\left\{ (-1)^r \int_{\real +i0}
\frac{e^{zt}}{t \prod_{j=1}^r (e^{\omega_j t}-1)} dt \right\}} \\
&=& \exp\left\{ (-1)^r 2 \pi i \sum_{k=1}^r \sum_{n=1}^\infty
\res_{t=\frac{2 \pi in}{\omega_k}}
\frac{e^{zt}}{t \prod_{j=1}^r (e^{\omega_j t}-1)} \right\} \\
&=& \exp\left\{ (-1)^r 2 \pi i \sum_{k=1}^r \sum_{n=1}^\infty
\frac{x_k^n}{2 \pi in \prod_{j=1, j \ne k}^r (q_{jk}^n-1)}
\right\} \\
&=& \prod_{k=1}^r \exp\left\{-\sum_{n=1}^\infty
\frac{x_k^n}{n \prod_{j=1, j \ne k}^r (1-q_{jk}^n)} \right\} \\
&=& \prod_{k=1}^r \exp\left\{
- \li_r (x_k ; \underline{q_k}) \right\} \\
&=& \prod_{k=1}^r (x_k ; \underline{q_k})^{(r-2)}_\infty.
\end{eqnarray*}
The last equality is due to Lemma \ref{3.6}.
Therefore we get (\ref{39}) from (\ref{32}) under the restriction
of $z$,
but we can conclude (\ref{39}) for any $z \in \comp$
except the poles by the analytic continuation.

Similarly in the domain
\[
\left\{ z \in \comp \Big|
0 < \re z < \re |\underline{\omega}| \ ,\ 
\im z <0 \ ,\ 
\im z < \im |\underline{\omega}| \ ,\ 
\im \frac{z}{\omega_k} <0 \ (\forall k)
\right\},
\]
we are allowed to add up the residues in $\{ \im t <0 \}$
and we get (\ref{40}) from (\ref{33}). \qed

\bigskip

Consequently it is clear that $S_r (z | \underline{\omega})$ is
meromorphic on $\comp$.
In the next section, Proposition \ref{3.8} will be connected
with the modular properties
of the multiple elliptic gamma functions
$G_r (z | \underline{\tau})$.

\begin{cor}\label{3.9}
If $\im \frac{\omega_1}{\omega_2} >0$, then
\begin{eqnarray*}
S_2 (z | \omega_1,\omega_2)
&=& \exp \left\{
+ \frac{\pi i}{2} B_{22} (z|\omega_1,\omega_2)
\right\}
\prod_{j=0}^{\infty}\frac{
1-e^{2 \pi i (z / \omega_2 +j \omega_1 / \omega_2)}
}{
1-e^{2 \pi i (z / \omega_1 -(j+1) \omega_2 / \omega_1)}
} \\
&=& \exp \left\{
- \frac{\pi i}{2} B_{22} (z|\omega_1,\omega_2)
\right\}
\prod_{j=0}^{\infty}\frac{
1-e^{2 \pi i ( - z / \omega_1 -j \omega_2 / \omega_1)}
}{
1-e^{2 \pi i ( - z / \omega_2 +(j+1) \omega_1 / \omega_2)}
}.
\end{eqnarray*}

If $\im \frac{\omega_1}{\omega_2},
\im \frac{\omega_1}{\omega_3},
\im \frac{\omega_2}{\omega_3} >0$, then
\begin{eqnarray*}
\lefteqn{S_3 (z | \omega_1,\omega_2,\omega_3)} \\
&=& \exp \left\{
- \frac{\pi i}{6} B_{33} (z|\omega_1,\omega_2,\omega_3)
\right\} \\
&& \times \prod_{j,k=0}^{\infty} \frac{
( 1-e^{2 \pi i
(z/\omega_1 -(j+1) \omega_2/\omega_1 -(k+1) \omega_3/\omega_1) } )
( 1-e^{2 \pi i
(z/\omega_3 +j \omega_1/\omega_3 +k \omega_2/\omega_3) } )
}{
1-e^{2 \pi i
(z/\omega_2 +j \omega_1/\omega_2 -(k+1) \omega_3/\omega_2) }
} \\
&=& \exp \left\{
+ \frac{\pi i}{6} B_{33} (z|\omega_1,\omega_2,\omega_3)
\right\} \\
&& \times \prod_{j,k=0}^{\infty} \frac{
( 1-e^{2 \pi i
(-z/\omega_1 -j \omega_2/\omega_1 -k \omega_3/\omega_1) } )
( 1-e^{2 \pi i
(-z/\omega_3 +(j+1) \omega_1/\omega_3 +(k+1) \omega_2/\omega_3) })
}{
1-e^{2 \pi i
(-z/\omega_2 +(j+1) \omega_1/\omega_2 -k \omega_3/\omega_2) }
}.
\end{eqnarray*}
\end{cor}

\bigskip

Similarly the equation
$ S_1 (z | \omega_1) = 2 \sin \frac{\pi z}{\omega_1} $
is confirmed by the residue formula.


\section{Modular properties of the multiple elliptic gamma
functions $G_r (z|\underline{\tau})$}

Now we show the main theorem of this paper.

\begin{thm}[Modular properties of $G_r (z| \underline{\tau})$]
\label{4.1}
If $r \ge 2, \im \frac{\omega_j}{\omega_k} \ne 0$,
then the multiple elliptic gamma function satisfies the identity
\[
\prod_{k=1}^r G_{r-2}
\left( \frac{z}{\omega_k} \bigg| \left(
\frac{\omega_1}{\omega_k}, \cdots,
\widecheck{\frac{\omega_k}{\omega_k}}, \cdots,
\frac{\omega_r}{\omega_k} \right) \right)
= \exp \left\{
- \frac{2 \pi i}{r!} B_{rr} (z|\underline{\omega})
\right\}.
\]
\end{thm}

\paragraph*{Proof}
Assume that $\omega_1,\cdots,\omega_r$ lie on the same
side of some line through the origin.
Then comparing the formulae (\ref{39}) and (\ref{40}) implies
\begin{eqnarray*}
\exp \left\{
- \frac{2 \pi i}{r!} B_{rr} (z|\underline{\omega}) \right\}
&=& \prod_{k=1}^{r}
\left\{ (x_k^{-1} ; \underline{q_k^{-1}})^{(r-2)}_\infty
 \right\}^{(-1)^{r-1}}
\left\{ (x_k ; \underline{q_k})^{(r-2)}_\infty \right\}^{(-1)^r}
 \\
&=& \prod_{k=1}^r G_{r-2}
\left( \frac{z}{\omega_k} \bigg| \left(
\frac{\omega_1}{\omega_k}, \cdots,
\widecheck{\frac{\omega_k}{\omega_k}}, \cdots,
\frac{\omega_r}{\omega_k} \right) \right),
\end{eqnarray*}
where we used the definition (\ref{13}) of
$G_r (z|\underline{\tau})$.

In general, it suffices to see the case when 
$\omega_2,\cdots,\omega_r$ lie on the same side and only
$\omega_1$ lies on the contour side.
In other words, $- \omega_1,\omega_2,\cdots,\omega_r$
lie on the same side.
Then apply (\ref{18}),(\ref{16}),(\ref{14}) and (\ref{23})
to complete the proof inductively.
In fact, it follows that
\begin{eqnarray*}
\lefteqn{
G_{r-2} \left( \frac{z}{\omega_1} \bigg| \left(
\frac{\omega_2}{\omega_1} , \cdots ,
\frac{\omega_r}{\omega_1} \right) \right)
\prod_{k=2}^r G_{r-2}
\left( \frac{z}{\omega_k} \bigg| \left(
\frac{\omega_1}{\omega_k} , \cdots ,
\widecheck{\frac{\omega_k}{\omega_k}} , \cdots ,
\frac{\omega_r}{\omega_k} \right) \right)
} \\
&=& G_{r-2} \left( \frac{z-\omega_1}{-\omega_1} \bigg| \left(
\frac{\omega_2}{-\omega_1} , \cdots ,
\frac{\omega_r}{-\omega_1} \right) \right)^{-1}
\prod_{k=2}^r G_{r-2}
\left( \frac{z-\omega_1}{\omega_k} \bigg| \left(
\frac{-\omega_1}{\omega_k} ,
\frac{\omega_2}{\omega_k} , \cdots ,
\widecheck{\frac{\omega_k}{\omega_k}} , \cdots ,
\frac{\omega_r}{\omega_k} \right) \right)^{-1} \\
&=& \exp \left\{
+ \frac{2 \pi i}{r!} B_{rr} (z-\omega_1|\underline{\omega}[1])
\right\} \\
&=& \exp \left\{
- \frac{2 \pi i}{r!} B_{rr} (z|\underline{\omega})
\right\}.
\qed
\end{eqnarray*}

\begin{thm}[Modular properties of $G_r (z| \underline{\tau})$]
\label{4.2}
If $\im\tau_j \ne 0$ and $\im\frac{\tau_j}{\tau_k} \ne 0$, then
\begin{eqnarray}
G_r (z| \underline{\tau})
&=& \exp \left\{
\frac{2 \pi i}{(r+2)!} B_{r+2,r+2} (z|(\underline{\tau},-1))
\right\} \nonumber \\
&& \times \prod_{k=0}^r G_r
\left( \frac{z}{\tau_k} \bigg| \left(
\frac{\tau_0}{\tau_k}, \cdots,
\widecheck{\frac{\tau_k}{\tau_k}}, \cdots,
\frac{\tau_r}{\tau_k},
-\frac{1}{\tau_k}
\right) \right) \label{41} \\
&=& \exp \left\{
- \frac{2 \pi i}{(r+2)!} B_{r+2,r+2} (z|(\underline{\tau},1))
\right\} \nonumber \\
&& \times \prod_{k=0}^r G_r
\left( -\frac{z}{\tau_k} \bigg| \left(
-\frac{\tau_0}{\tau_k}, \cdots,
-\widecheck{\frac{\tau_k}{\tau_k}}, \cdots,
-\frac{\tau_r}{\tau_k},
-\frac{1}{\tau_k}
\right) \right). \label{42}
\end{eqnarray}
\end{thm}

\paragraph*{Proof}
In Theorem \ref{4.1}, we replace $r$ by $r+2$ to get
\begin{eqnarray*}
\lefteqn{
\prod_{k=0}^{r+1} G_r
\left( \frac{z}{\omega_k} \bigg| \left(
\frac{\omega_0}{\omega_k}, \cdots,
\widecheck{\frac{\omega_k}{\omega_k}}, \cdots,
\frac{\omega_{r+1}}{\omega_k} \right) \right)
} \\
&=& \exp \left\{ - \frac{2 \pi i}{(r+2)!}
B_{r+2,r+2} (z|( \omega_0, \cdots, \omega_{r+1})) \right\}.
\end{eqnarray*}
When $\omega_0=\tau_0, \cdots, \omega_r=\tau_r, \omega_{r+1}=-1$,
this becomes
\begin{eqnarray*}
\lefteqn{
G_r (-z|-\underline{\tau})
\prod_{k=0}^r G_r
\left( \frac{z}{\tau_k} \bigg| \left(
\frac{\tau_0}{\tau_k}, \cdots,
\widecheck{\frac{\tau_k}{\tau_k}\ }, \cdots,
\frac{\tau_r}{\tau_k},
-\frac{1}{\tau_k}
\right) \right)
} \\
&=& \exp \left\{ - \frac{2 \pi i}{(r+2)!}
B_{r+2,r+2} (z|(\underline{\tau},-1)) \right\}.
\end{eqnarray*}
The first formula follows from (\ref{18}).

Second formula of the theorem is similar to first one by letting
$\omega_{r+1} = +1$.
\qed

\bigskip

The equivalence between (\ref{41}) and (\ref{42}) follows from
(\ref{29}),(\ref{18}),(\ref{19}) and Theorem \ref{4.1}.
Of course, this theorem includes Jacobi's result (\ref{3})
and Felder and Varchenko's result (\ref{1}).

\begin{cor}
If $\im\tau >0$, then
\begin{eqnarray*}
\theta_0 (z,\tau)
&=& \exp \{ \pi i B_{22} (z | \tau,-1) \}
\theta_0 \left( \frac{z}{\tau},-\frac{1}{\tau} \right) \\
&=& \exp \{ - \pi i B_{22} (z | \tau,1) \}
\theta_0 \left( - \frac{z}{\tau},-\frac{1}{\tau} \right).
\end{eqnarray*}

If $\im\tau, \im\sigma, \im\frac{\tau}{\sigma} >0$, then
\begin{eqnarray*}
\Gamma (z,\tau,\sigma)
&=& \exp \left\{
\frac{\pi i}{3} B_{33} (z | \tau,\sigma,-1) \right\}
\frac{
\Gamma \left( \frac{z}{\sigma},
\frac{\tau}{\sigma},-\frac{1}{\sigma} \right)
}{
\Gamma \left( \frac{z-\sigma}{\tau},
-\frac{\sigma}{\tau},-\frac{1}{\tau} \right)
} \\
&=& \exp \left\{
- \frac{\pi i}{3} B_{33} (z | \tau,\sigma,1) \right\}
\frac{
\Gamma \left( - \frac{z}{\tau},
-\frac{\sigma}{\tau},-\frac{1}{\tau} \right)
}{
\Gamma \left( \frac{\tau -z}{\sigma},
\frac{\tau}{\sigma},-\frac{1}{\sigma} \right)
}.
\end{eqnarray*}

If $\im\tau_0, \im\tau_1, \im\tau_2, \im\frac{\tau_0}{\tau_1},
\im\frac{\tau_0}{\tau_2}, \im\frac{\tau_1}{\tau_2} >0$, then
\begin{eqnarray*}
G_2 (z|\tau_0,\tau_1,\tau_2)
&=& \exp \left\{
\frac{\pi i}{12} B_{44} (z | \tau_0,\tau_1,\tau_2,-1) \right\} \\
&& \times \frac{
G_2 \left( \dfrac{z-\tau_1 -\tau_2}{\tau_0} \bigg|
- \frac{\tau_1}{\tau_0}, - \frac{\tau_2}{\tau_0},
- \frac{1}{\tau_0} \right)
G_2 \left( \dfrac{z}{\tau_2} \bigg| \frac{\tau_0}{\tau_2},
\frac{\tau_1}{\tau_2},-\frac{1}{\tau_2} \right)
}{
G_2 \left( \dfrac{z-\tau_2}{\tau_1} \bigg|
\frac{\tau_0}{\tau_1}, -\frac{\tau_2}{\tau_1},
-\frac{1}{\tau_1} \right)
} \\
&=& \exp \left\{
- \frac{\pi i}{12} B_{44} (z | \tau_0,\tau_1,\tau_2,1) \right\} \\
&& \times \frac{
G_2 \left( -\dfrac{z}{\tau_0} \bigg|
-\frac{\tau_1}{\tau_0},-\frac{\tau_2}{\tau_0},
-\frac{1}{\tau_0} \right)
G_2 \left( \dfrac{\tau_0 + \tau_1 -z}{\tau_2} \bigg|
\frac{\tau_0}{\tau_2},\frac{\tau_1}{\tau_2},
-\frac{1}{\tau_2} \right)
}{
G_2 \left( \dfrac{\tau_0 -z}{\tau_1} \bigg|
\frac{\tau_0}{\tau_1},-\frac{\tau_2}{\tau_1},
-\frac{1}{\tau_1} \right)
}.
\end{eqnarray*}
The parameters of $\theta_0, \Gamma, G_2$ which appeared above
are all in the upper half plane.
\end{cor}

\paragraph*{Proof}
To prove the first identity of $G_2$,
we take $r=2$ in Theorem \ref{4.2}, and get
\begin{eqnarray*}
G_2 (z | \tau_0,\tau_1,\tau_2)
&=& \exp \left\{
\frac{\pi i}{12} B_{44} (z | \tau_0,\tau_1,\tau_2,-1) \right\}
G_2 \left( \frac{z}{\tau_0} \bigg|
\frac{\tau_1}{\tau_0},\frac{\tau_2}{\tau_0},
-\frac{1}{\tau_0} \right) \\
&& \times
G_2 \left( \frac{z}{\tau_1} \bigg|
\frac{\tau_0}{\tau_1},\frac{\tau_2}{\tau_1},
-\frac{1}{\tau_1} \right)
G_2 \left( \frac{z}{\tau_2} \bigg|
\frac{\tau_0}{\tau_2},\frac{\tau_1}{\tau_2},
-\frac{1}{\tau_2} \right).
\end{eqnarray*}
We recall that the functional equation (\ref{16}) implies
\begin{eqnarray*}
G_2 \left( \frac{z}{\tau_0} \bigg|
\frac{\tau_1}{\tau_0},\frac{\tau_2}{\tau_0},
-\frac{1}{\tau_0} \right)
&=& G_2 \left( \frac{z-\tau_1 -\tau_2}{\tau_0} \bigg|
- \frac{\tau_1}{\tau_0}, - \frac{\tau_2}{\tau_0},
- \frac{1}{\tau_0} \right), \\
G_2 \left( \frac{z}{\tau_1} \bigg|
\frac{\tau_0}{\tau_1},\frac{\tau_2}{\tau_1},
-\frac{1}{\tau_1} \right)
&=& G_2 \left( \dfrac{z-\tau_2}{\tau_1} \bigg|
\frac{\tau_0}{\tau_1}, -\frac{\tau_2}{\tau_1},
-\frac{1}{\tau_1} \right)^{-1}.
\end{eqnarray*}
All required formulae are obtained in a similar way.
\qed


\section{Representation of $G_r (z|\underline{\tau})$
by the integral or $S_{r+1} (z|\underline{\omega})$}

In this section, we start with the following proposition,
which is verified by evaluating the residues
inside the contour $C_1$.

\begin{prop}\label{5.2}
When $\im z >0, \im\tau_j \ne 0 \ (\forall j)$, then we have
\[
\li_{r+2} (x;\underline{q})
= - \int_{C_1} \frac{ e^{2 \pi izt} }{
t (1- e^{2 \pi it}) \prod_{j=0}^r (1- e^{2 \pi i \tau_j t})
} dt,
\]
where the contour $C_1$ is as follows.

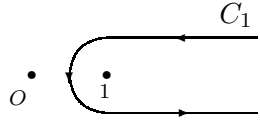
\begin{center}
\setlength{\unitlength}{1mm}
\begin{picture}(40,20)(-10,-10)
\put(0,0){\circle*{1}}
\put(-3,-4){$\scriptstyle O$}
\put(10,0){\circle*{1}}
\put(9,-3){$\scriptstyle 1$}
\put(10,5){\line(1,0){20}}
\put(10,-5){\line(1,0){20}}
\bezier{100}(5,0)(5.3,4.7)(10,5)
\bezier{100}(5,0)(5.3,-4.7)(10,-5)
\put(20,5){\vector(-1,0){1}}
\put(20,-5){\vector(1,0){1}}
\put(5,0){\vector(0,-1){1}}
\put(25,7){$C_1$}
\end{picture}
\end{center}
\end{prop}

The definition of $G_r (z|\underline{\tau})$,
Lemma \ref{3.6} and the above proposition lead us
to the following representations.
\begin{thm}[Integral representations of
$G_r (z|\underline{\tau})$]
When $\im\tau_j >0 \ (\forall j),
\ 0< \im z < \im |\underline{\tau}| $, then
\begin{eqnarray*}
\lefteqn{G_r (z|\underline{\tau})} \\
&=& \exp\left\{ \int_{C_1} \frac{
e^{2 \pi izt} +(-1)^r e^{2 \pi i (|\underline{\tau}|+1-z) t}
}{
t (e^{2 \pi it} -1) \prod_{j=0}^r (e^{2 \pi i \tau_j t} -1)
} dt \right\} \\
&=& \exp\left\{
- \frac{2 \pi i}{(r+2)!} B_{r+2,r+2} (z|(\underline{\tau},1))
+ \int_{\real +i \epsilon} \frac{
-e^{2 \pi izt} +(-1)^{r+1} e^{2 \pi i (|\underline{\tau}| +1-z) t}
}{
t (e^{2 \pi it} -1) \prod_{j=0}^r (e^{2 \pi i \tau_j t} -1)
} dt \right\} \\
&=& \exp\left\{
\frac{2 \pi i}{(r+2)!} B_{r+2,r+2} (z|(\underline{\tau},-1))
+ \int_{\real +i \epsilon} \frac{
e^{2 \pi izt} +(-1)^r e^{2 \pi i (|\underline{\tau}| -1-z) t}
}{
t (e^{-2 \pi it} -1) \prod_{j=0}^r (e^{2 \pi i \tau_j t} -1)
} dt \right\}. \label{46}
\end{eqnarray*}
where we take a small $\epsilon >0$ and
the contours as following figures.

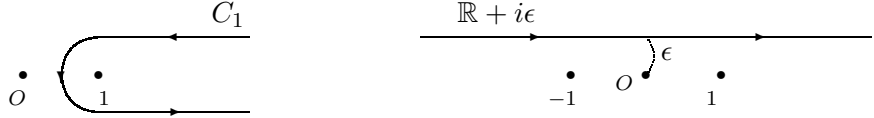
\begin{center}
\setlength{\unitlength}{1mm}
\begin{picture}(40,20)(-10,-10)
\put(0,0){\circle*{1}}
\put(-2,-4){$\scriptstyle O$}
\put(10,0){\circle*{1}}
\put(10,-4){$\scriptstyle 1$}
\put(10,5){\line(1,0){20}}
\put(10,-5){\line(1,0){20}}
\bezier{100}(5,0)(5.3,4.7)(10,5)
\bezier{100}(5,0)(5.3,-4.7)(10,-5)
\put(20,5){\vector(-1,0){1}}
\put(20,-5){\vector(1,0){1}}
\put(5,0){\vector(0,-1){1}}
\put(25,7){$C_1$}
\end{picture}
\hspace{20mm}
\begin{picture}(60,20)(-30,-10)
\put(0,0){\circle*{1}}
\put(-4,-2){$\scriptstyle O$}
\put(10,0){\circle*{1}}
\put(8,-4){$\scriptstyle 1$}
\put(-10,0){\circle*{1}}
\put(-13,-4){$\scriptstyle -1$}
\put(-25,7){$\real +i \epsilon$}
\put(-30,5){\line(1,0){60}}
\put(-15,5){\vector(1,0){1}}
\put(15,5){\vector(1,0){1}}
\bezier{15}(0,5)(2,2.5)(0,0)
\put(2,2){$\epsilon$}
\end{picture}
\end{center}
\end{thm}

We have immediate results by setting $r=0$ or $1$.

\begin{cor}
If $0 < \im z < \im \tau$, we have
\begin{eqnarray*}
\theta_0 (z,\tau)
&=& \exp\left\{ - \int_{C_1} \frac{
\cos ( \pi (2z -\tau -1) t)
}{
2t \sin (\pi t) \sin (\pi \tau t)
} dt \right\} \\
&=& \exp\left\{ - \pi i B_{22} (z|\tau,1)
+ \int_{\real +i \epsilon} \frac{
\cos ( \pi (2z -\tau -1) t)
}{
2t \sin (\pi t) \sin (\pi \tau t)
} dt \right\} \\
&=& \exp\left\{ \pi i B_{22} (z|\tau,-1)
+ \int_{\real +i \epsilon} \frac{
\cos ( \pi (2z -\tau +1) t)
}{
2t \sin (\pi t) \sin (\pi \tau t)
} dt \right\}.
\end{eqnarray*}

If $0 < \im z < \im (\tau+\sigma)$, we have
\begin{eqnarray*}
\Gamma (z,\tau,\sigma)
&=& \exp\left\{ - \int_{C_1} \frac{
\sin ( \pi (2z -\tau -\sigma -1) t)
}{
4t \sin (\pi t) \sin (\pi \tau t) \sin (\pi \sigma t)
} dt \right\} \label{47} \\
&=& \exp\left\{ - \frac{\pi i}{3} B_{33} (z|\tau,\sigma,1)
+ \int_{\real +i \epsilon} \frac{
\sin ( \pi (2z -\tau -\sigma -1) t)
}{
4t \sin (\pi t) \sin (\pi \tau t) \sin (\pi \sigma t)
} dt \right\} \\
&=& \exp\left\{ \frac{\pi i}{3} B_{33} (z|\tau,\sigma,-1)
+ \int_{\real +i \epsilon} \frac{
\sin ( \pi (2z -\tau -\sigma +1) t)
}{
4t \sin (\pi t) \sin (\pi \tau t) \sin (\pi \sigma t)
} dt \right\}.
\end{eqnarray*}
\end{cor}

The corollary implies the following formulae
by collecting the residues inside the contour $C_1$.
These were originally obtained
in a different way \cite{FV}.
Generally $G_r (z|\underline{\tau})$ have
similar representations \cite{N1}.

\begin{cor}[summation formula]
If $0 < \im z < \im \tau$, then
\[
\theta_0 (z,\tau)
= \exp \left( -i \sum_{j=1}^{\infty}
\frac{\cos (\pi j(2z - \tau))}
{j \sin (\pi j \tau)} \right).
\]

If $0 < \im z < \im (\tau+\sigma)$, then
\[
\Gamma (z,\tau,\sigma)
= \exp \left( -\frac{i}{2} \sum_{j=1}^{\infty}
\frac{\sin (\pi j(2z - \tau - \sigma))}
{j \sin (\pi j \tau) \sin (\pi j \sigma)} \right).
\]
\end{cor}

At the end of this paper, we give the following theorem
which comes from the integral representations of
$G_r (z| \underline{\tau})$ and $S_r (z| \underline{\omega})$.
We also use (\ref{20}) and Proposition \ref{3.1}(i) for the proof.

\begin{thm}[Representations of $G_r$ by the infinite product
of $S_{r+1}$]\label{5.6}
We assume $\im\tau_j >0 \ (\forall j)$,
$0< \im z < \im |\underline{\tau}| $, then
\begin{eqnarray*}
G_r (z| \underline{\tau})
&=& \exp \left\{
\frac{2 \pi i}{(r+2)!} B_{r+2,r+2} (z|(\underline{\tau},-1))
\right\} \nonumber\\
&& \times \prod_{k=0}^{\infty} \frac{
S_{r+1} (z+k+1|\underline{\tau})^{(-1)^r}
S_{r+1} (z-k |\underline{\tau})^{(-1)^r}
}{
\exp \left\{ \frac{\pi i}{(r+1)!} (
 B_{r+1,r+1} (z+k+1|\underline{\tau})
- B_{r+1,r+1} (z-k |\underline{\tau})
) \right\}
} \label{50} \\
&=& \exp \left\{
- \frac{2 \pi i}{(r+2)!} B_{r+2,r+2} (z|(\underline{\tau},1))
\right\} \nonumber \\
&& \times \prod_{k=0}^{\infty} \frac{
S_{r+1} (z+k |\underline{\tau})^{(-1)^r}
S_{r+1} (z-k-1|\underline{\tau})^{(-1)^r}
}{
\exp \left\{\frac{\pi i}{(r+1)!} (
 B_{r+1,r+1} (z+k |\underline{\tau})
- B_{r+1,r+1} (z-k-1|\underline{\tau})
) \right\}
}.
\end{eqnarray*}
\end{thm}

\bigskip

Now we denote
\[
\psi_2 (z) = \exp \left( 2 \pi i
\int_{-i \infty}^{z} \frac{t-1}{e^{2 \pi it}-1} \ dt \right),
\]
which satisfies
$\psi_2 (1) = \exp \frac{\pi i}{12}$ and
$\psi_2 (z) \psi_2 (2-z)
= \exp \left( - \pi i B_{22} (z|1,1) \right)$.
Refering to Kurokawa \cite{K1,K2}, $S_2 (z | 1,1)$ obeys
the identity
\[
S_2 (z | 1,1)
= \exp \left( - \int_1^z \pi (t-1) \cot (\pi t) \ dt \right).
\]
Hence they imply
\begin{eqnarray*}
S_2 (z | 1,1)
&=& \psi_2 (z)^{-1}
\exp \left\{ - \frac{\pi i}{2} B_{22} (z | 1,1) \right\} \\
&=& \psi_2 (2-z)
\exp \left\{ \frac{\pi i}{2} B_{22} (z | 1,1) \right\}.
\end{eqnarray*}

Substituting $r=1,\  \tau = \tau_0 = \tau_1$ in Theorem \ref{5.6},
we can get the formula
\[
\Gamma (z,\tau,\tau)
= \exp\left\{ \frac{\pi i}{3} B_{33} (z | \tau,\tau,-1) \right\}
\prod_{k=0}^{\infty} \frac{
\psi_2 \left( \frac{z+k+1}{\tau} \right)
}{
\psi_2 \left( 2-\frac{z-k}{\tau} \right)
}.
\]
This is a shifted version of the formula which appeared in
{\it Theorem 5.2} of paper \cite{FV}.


\paragraph*{Acknowledgment}
The author thanks Prof. Kimio Ueno
for discussions and support.
Thanks are also due to Dr. Michitomo Nishizawa.


\end{document}